\documentclass[a4paper,11pt,twoside,reqno]{amsart}

\usepackage[utf8]{inputenc}
\usepackage[plainpages=false,pdfpagelabels=true]{hyperref}
\usepackage{amsmath}
\usepackage{amssymb,amsthm}
\usepackage[margin=1in]{geometry}
\usepackage{slashed}

\newtheorem{Satz}{Theorem}[section]
\newtheorem{Prop}[Satz]{Proposition}
\newtheorem{Lem}[Satz]{Lemma}

\newtheorem{Cor}[Satz]{Corollary}

\theoremstyle{definition}

\newtheorem{Bem}[Satz]{Remark}

\newcommand{\tr}{\operatorname{Tr}}

\newcommand{\dv}{\text{ }dV}
\newcommand{\s}{{\mathbb{S}}}
\allowdisplaybreaks[1]

\renewcommand{\epsilon}{\varepsilon}

\newcommand{\R}{\ensuremath{\mathbb{R}}}
\newcommand{\N}{\ensuremath{\mathbb{N}}}

\numberwithin{equation}{section}

\title{On polyharmonic helices in space forms}
\author{Volker Branding}
\date{\today}
\address{University of Vienna, Faculty of Mathematics\\
Oskar-Morgenstern-Platz 1, 1090 Vienna, Austria\\}
\email{volker.branding@univie.ac.at}
\thanks{The author gratefully acknowledges the support of the Austrian Science Fund (FWF) via the START-project Y963-N35 of Michael Eichmair and 
the project P34853-N35 “Geometric Analysis of Biwave Maps".}

\begin{document}

\begin{abstract}
In this article we study polyharmonic curves of constant curvature
where we mostly focus on the case of curves on the sphere.
We classify polyharmonic curves
of constant curvature in three-dimensional space forms and
derive an explicit family of polyharmonic curves
on the sphere. 
Our results give new insights into the geometric structure
of higher order variational problems.
\end{abstract} 

\maketitle

\section{Introduction and Results}
Geometric variational problems are a vibrant field of research with many links to analysis, geometry and theoretical physics.
One of their most prominent representative
are \emph{harmonic maps} which are a
subject on which a significant number of results could be achieved by many mathematicians.

Let us briefly recall the mathematical setup employed in the theory of harmonic maps.
We consider a smooth map \(\phi\) between two Riemannian manifolds \((M,g)\) and \((N,h)\)
and its associated energy 
\begin{align}
\label{energy-map}
E_1(\phi)=E(\phi):=\int_M|d\phi|^2\dv.
\end{align}
The critical points of \eqref{energy-map} are characterized by the vanishing of the so-called
\emph{tension field} which is defined as follows
\begin{align}
\label{tension-field}
\tau(\phi):=\tr_g\bar\nabla d\phi \in\Gamma(\phi^\ast TN),
\end{align}
where \(\bar\nabla\) represents the connection on \(\phi^\ast TN\).
The solutions of \(\tau(\phi)=0\) are called \emph{harmonic maps}.
The harmonic map equation is a second order semilinear elliptic partial differential equation for which
many powerful tools, such as the maximum principle,
can be used to investigate existence and qualitative behavior of a given solution.

There are various variants in the mathematics literature that extend the classical energy functional
\eqref{energy-map} taking into account higher derivatives of the map \(\phi\).

One possibility is given by the following energy functional, usually referred to as \emph{\(r\)-energy},
where we need to distinguish between the even and the odd case
\begin{align}
\label{energy-r-harmonic}
E_{2s}(\phi):=&\int_M|\underbrace{(d^\ast d)\ldots (d^\ast d)}_{s-\textrm{times}}\phi|^2\dv,\qquad s=1,\ldots \\
\nonumber E_{2s+1}(\phi):=&\int_M|d\underbrace{(d^\ast d)\ldots (d^\ast d)}_{s-\textrm{times}}\phi|^2\dv,\qquad s=1,\ldots.
\end{align}
The critical points of \eqref{energy-r-harmonic} are called \emph{polyharmonic maps of order \(r\)}
or simply \emph{\(r\)-harmonic maps}
and were explicitly determined by Wang \cite{wang}
and in addition by Maeta \cite{MR3007953}, 
who also calculated the second variation of \eqref{energy-r-harmonic}.

Polyharmonic maps are a semilinear elliptic equation of order \(2r\) making
it a challenging problem to find explicit solutions of the latter.
For the explicit structure of \(r\)-harmonic maps we refer to \cite{MR3007953}.

Let us mention several recent results obtained for the critical points of \eqref{energy-r-harmonic}.
Explicit solutions of the Euler-Lagrange equations of \eqref{energy-r-harmonic} were constructed in
\cite{MR3711937,MR3790367}. Hypersurfaces in space forms characterized by critical points of \eqref{energy-r-harmonic}
have been investigated in \cite{bmor2021poly,MR4462636}.
For classification results of the critical points of \eqref{energy-r-harmonic} we refer to \cite{MR4184658}
and various unique continuation properties were studied in \cite{bmor2021}.
The stress-energy tensor of \eqref{energy-r-harmonic} has been systematically analyzed in \cite{MR4007262}.

Another possibility of obtaining a higher order energy functional for maps between Riemannian manifolds
was suggested by Eells and Sampson in 1964 \cite{MR172310} by considering 

\begin{align}
\label{energy-es}
E_r^{ES}(\phi):=\int_M|(d+d^\ast)^r\phi|^2\dv,\qquad r=1,\ldots.
\end{align}

The critical points of \eqref{energy-es} are referred to as \emph{ES-r-harmonic maps}.
In general, the critical points of \eqref{energy-r-harmonic} and \eqref{energy-es} will be given by a different
set of equations.
However, if \(r=1,2,3\) or in the case of a one-dimensional domain both functionals \eqref{energy-r-harmonic}
and \eqref{energy-es} coincide.
For an extensive analysis of \eqref{energy-es} and a discussion of the differences between \eqref{energy-r-harmonic}
and \eqref{energy-es} we refer to the recent articles \cite{MR4106647,MR4293944,bmor2021,MR4462636}.

A drawback of the energy functionals \eqref{energy-r-harmonic} and \eqref{energy-es} is the fact that they are not coercive
which makes it a challenging problem to prove the existence of critical points.
However, there exists another higher order energy functional which overcomes this problem.
In this case it is assumed that the target manifold \(N\) is realized as a submanifold of some \(\R^q\).
For example, in the case of a spherical target we can consider smooth maps
\(u\colon M\to\s^n\subset\R^{n+1}\) and in this case the energy functional for
\emph{extrinsic polyharmonic maps} is given by 
\begin{align}
\label{energy-extrinsic}
E^{ext}_{2s}(u):=&\int_M|\Delta^su|^2\dv,\qquad s=1,\ldots, \\
\nonumber E^{ext}_{2s+1}(u):=&\int_M|\nabla\Delta^su|^2\dv,\qquad s=1,\ldots.
\end{align}
For a recent analysis of \eqref{energy-extrinsic} and related background material we refer to  
\cite{MR4436204}.

In this article we will consider the aforementioned energy functionals \eqref{energy-r-harmonic} and \eqref{energy-es}
in the case of a one-dimensional domain and 
mostly choose a spherical target.

The \(r\)-energy for a curve \(\gamma\colon I\subset\R \to N\) (corresponding to both \eqref{energy-r-harmonic} and \eqref{energy-es} if \(\dim M=1\)) is given by  
\begin{align}
\label{energy-poly-curve}
E_r(\gamma):=\int_I|\nabla_T^{r-1}T|^2ds,
\end{align}
where \(T=\gamma'\) is the tangent vector of \(\gamma\)
and \(s\) represents the parameter of the curve \(\gamma\).

The critical points of \eqref{energy-poly-curve} were calculated in \cite{MR3007953,wang} and are
characterized by the equation
\begin{align}
\label{el-poly}
\tau_r(\gamma)=\nabla^{2r-1}_TT+\sum_{l=0}^{r-2}(-1)^lR^N(\nabla_T^{2r-3-l}T,\nabla_T^lT)T=0
\end{align}
with \(R^N\) being the curvature tensor of the manifold \(N\).

Solutions of \eqref{el-poly} are called \emph{polyharmonic curves of order \(r\)} or shortly \emph{r-harmonic curves}.
In the case of \(r=1\) the energy \eqref{energy-poly-curve} reduces to the usual energy of a curve
whose critical points are \emph{geodesics}. Clearly, every geodesic is a solution of the equation for \(r\)-harmonic
curves \eqref{el-poly}, hence we are interested in finding non-geodesic solutions of \eqref{el-poly}
which we will call \emph{proper} \(r\)-harmonic curves.

Let us give an overview on the results for \(r\)-harmonic curves
which have already been established in the mathematics literature.
Much research has been performed for \(r=2\) in which case solutions of \eqref{el-poly}
are called \emph{biharmonic curves}. 
Biharmonic curves on surfaces were investigated in \cite{MR1884940},
for biharmonic curves on the Euclidean sphere we refer to
\cite{MR2187367,MR1863283,MR1919374}.
The Euler-Lagrange method for biharmonic curves was introduced in
\cite{MR2274737} and further developed in \cite{MR3045700}.
Making use of an algebraic approach proper biharmonic curves on quadrics were studied in \cite{MR3292682},
explicit formulas for proper biharmonic curves in Sasakian space forms have been obtained in \cite{MR2485475}.
For an extended summary on biharmonic curves we refer to \cite[Section 4]{MR3098705}.

If \(r=3\) the solutions of \eqref{el-poly} are called \emph{triharmonic curves}.
Some general results on triharmonic maps have been obtained in \cite{MR3403738},
triharmonic curves on surfaces and three dimensional spaces with constant curvature were studied in \cite{MR3007953}.
Recently, triharmonic curves on three-dimensional homogeneous spaces \cite{montri}
and on \(f\)-Kenmotsu Manifolds \cite{boz2021}
have been investigated.
While biharmonic curves always need to have constant geodesic curvature 
it was shown in \cite[Theorem 1.1]{montri} that there exist proper triharmonic curves whose
geodesic curvature is non-constant.

Not many results have been established for \(r\)-harmonic curves
and \(r\geq 3\) as this comes
with many challenging technical difficulties. As already mentioned above a geodesic
is always a solution of the equation for \(r\)-harmonic curves \eqref{el-poly}.
However, it was shown by Maeta \cite{MR3007953} that biharmonic curves are not triharmonic in general
and we thus expect that \(r\)-harmonic curves will be different from \(q\)-harmonic curves, whenever \(q\neq r\).
For the investigation of \(r\)-harmonic curves in Euclidean spaces we refer to \cite[Section 3]{MR2869168}.

It should also be mentioned that it is possible to produce \(r\)-harmonic curves
from geodesics by reparametrizing the domain using a diffeomorphism.
To this end assume that \(\varphi(t)\) is a geodesic,
we reparametrize it using a diffeomorphism \(\mu(s)\)
and consider \(\psi(s):=(\varphi\circ\mu)(s)\).
Then, a direct calculation shows that \(\psi(s)\) is a proper
\(r\)-harmonic curve if and only if \(\mu(s)\) is a polynomial of order \(r'\), where \(2\leq r'\leq 2r-1\).
For more details, see \cite[Section 2.1]{MR4106647}.

Throughout this article we will use the word helix in order
to denote a curve whose geodesic curvature \(k\) and whose torsion \(\tau\) are  constant.

The first main contribution of this article is the following theorem 
characterizing \(r\)-harmonic curves on three-dimensional space forms.
Note that in the following theorem the symbol \(\tau\) represents
the torsion of the curve \(\gamma\) and does not represent the tension field
as defined in \eqref{tension-field}.

\begin{Satz}
\label{theorem-relation-kt-3d}
Let \(\gamma\colon I\to N\) be an \(r\)-harmonic curve parametrized by arclength 
where \(N\) is a three-dimensional space form
with constant curvature \(K\).
Furthermore, we assume that the geodesic curvature \(k\)
and the torsion \(\tau\) of the curve \(\gamma\) are constant.
Then the following relation holds
\begin{align}
\label{relation-kt-3d}
(k^2+\tau^2)^2=K\big((r-1)k^2+\tau^2\big).
\end{align}
\end{Satz}
This result extends a corresponding result for triharmonic curves \cite[Theorem 4.4]{montri} to arbitrary values of \(r\). Moreover, a similar statement for \(r\)-harmonic curves on surfaces was given in \cite[Proposition 5.6, Proposition 5.7]{MR3007953}.

A direct consequence of Theorem \ref{theorem-relation-kt-3d} is the following
\begin{Cor}
Let \(\gamma\colon I\to N\) be an \(r\)-harmonic curve parametrized by arclength 
where \(N\) is a three-dimensional space form
with constant non-positive curvature \(K\leq 0\).
Furthermore, we assume that the geodesic curvature \(k\)
and the torsion \(\tau\) of the curve \(\gamma\) are constant. Then \(\gamma\) is a geodesic.
\end{Cor}

\begin{Bem}
\begin{enumerate}
 \item The last statement supports the fact that for a target with negative curvature
polyharmonic maps necessarily have to be harmonic while in the case of a spherical
target there may be additional solutions besides harmonic maps.
\item If we consider an \(r\)-harmonic helix in the sphere \(\s^n\) 
with \(n\) sufficiently large, one can show by rewriting the 
Euler-Lagrange equation \eqref{el-poly} in terms of its Frenet-frame
that \(r\)-harmonic helices actually lie in \(\s^{2r-1}\).
Here, we use the terminology helix to denote a curve of which all geodesic
curvatures in its Frenet-frame are constant.
As \(\s^{2r-1}\) is totally geodesic in \(\s^n\) for \(n\) sufficiently large,
this supports the fact that the most general \(r\)-harmonic helix of a sphere
lies in \(\s^{2r-1}\).
\end{enumerate}
\end{Bem}

Using the explicit parametrization of helices in \(\s^3\) 
we can get the following application of Theorem \ref{theorem-relation-kt-3d}:

\begin{Cor}
\label{cor-helices}
For any \(r\geq 3\) there exists a family of \(r\)-harmonic curves \(\gamma\colon I\to\s^3\) with constant positive geodesic curvature \(k\) and constant positive torsion \(\tau\).
\end{Cor}

The second main result of this article is the following existence result for \(r\)-harmonic curves on the sphere.
\begin{Satz}
\label{theorem-r-harmonic-a}
The curve \(\gamma\colon I\to\s^n\) given by
\begin{align}
\label{solution-kharmonic-a}
\gamma(s)=\cos(\sqrt{r}s)e_1+\sin(\sqrt{r}s)e_2+e_3,
\end{align}
where \(e_i,i=1,2,3\) are mutually perpendicular and satisfy \(|e_1|^2=|e_2|^2=\frac{1}{r},|e_3|^2=\frac{r-1}{r}\)
is a proper \(r\)-harmonic curve which is parametrized by arclength.
\end{Satz}

\begin{Bem}
\hspace{1cm}
\label{remark-general-solution}
\begin{enumerate}
\item 
The curves \eqref{solution-kharmonic-a} are planar and have constant geodesic curvature 
\(k^2=r-1\) which is precisely what we expect from \eqref{relation-kt-3d} on the two-dimensional sphere \(\s^2\).

 \item The family of \(r\)-harmonic curves obtained in Theorem \ref{theorem-r-harmonic-a} is unstable in the sense that for such curves the second variation of the \(r\)-energy is negative.
 A detailed investigation of the stability of \(r\)-harmonic curves 
 \(\gamma\colon I\to\s^2\)
for \(r=2\) was carried out in \cite{MR4216418} and for \(r=2,3,4\) 
in \cite[Section 5]{MR4106647}.

\item 
Since the classification of \(r\)-harmonic curves is up to isometries
we can assume that \(k>0\) and \(\tau\geq 0\).
Let us denote \(x=k^2\) and \(y=\tau^2\) then \eqref{relation-kt-3d} with \(K=1\),
that is the case of a spherical target, will give us the equation of a conic and we are interested
in the intersection of the conic with the first quadrant.
Actually, when \(r=2\), the conic consists of two parallel lines 
given by \(x+y=0\) and \(x+y=1\) and the biharmonic curves
can be indexed by the segment \((A,B]\), where \(A(0,1)\) and \(B(1,0)\).
In the case that \(r\geq 3\) the conic is a parabola and its intersection with the first
quadrant is a non-empty and connected segment of a parabola.

\item It is well-known that a helix in \(\s^3\) is a geodesic on the generalized
Clifford torus \(T=\s^1(\cos(\alpha))\times\s^1(\sin(\alpha)),\alpha\in(0,\pi/2)\). In the case that \(\gamma\) is a proper biharmonic curve in \(\s^3\), then it
is a geodesic on the Clifford torus \(T=\s^1(\frac{1}{\sqrt{2}})\times\s^1(\frac{1}{\sqrt{2}})\subset\s^3\), see for example \cite[Theorem 2.2]{MR2187367}.
Note that this torus is special as it is minimal in \(\s^3\).
However, for an \(r\)-harmonic helix in \(\s^3\) with \(r\geq 3\) arbitrary, but fixed, the situation is different.
It lies on a torus \(T=\s^1(\cos(\alpha))\times\s^1(\sin(\alpha))\),
which is never the Clifford torus,
and which is different for each member of the family of \(r\)-harmonic helices
described in Corollary \ref{cor-helices} while in the case of biharmonic curves
the torus is always the Clifford torus. 
\item The last two items strongly suggest that biharmonic curves
are special among all \(r\)-harmonic curves.
\end{enumerate}
\end{Bem}

The statement of Theorem \ref{theorem-r-harmonic-a}
is known in the mathematics literature \cite[Proposition 4.4]{MR1919374}
in the biharmonic case (\(r=2\)).

Similar results for the inclusion map \(\iota:\s^n\hookrightarrow\s^{n+1}\) were already
obtained.
It was shown in \cite[Theorem 1.1]{MR3711937} that the canonical inclusion 
\(\iota:\s^n(R)\hookrightarrow\s^{n+1}\) is a proper critical point of \eqref{energy-r-harmonic},
which corresponds to a proper \(r\)-harmonic submanifold \(\s^n(R)\) of \(\s^{n+1}\)
if and only if \(R=\frac{1}{\sqrt{r}}\). This result was later generalized 
to critical points of \eqref{energy-es} in \cite[Theorem 1.1]{MR4106647}.

Throughout this article we will use the following notation.
By \(s\) we will denote the parameter of the curve \(\gamma\),
the first, second and third derivative of \(\gamma\) will be written as \(T:=\gamma',\gamma''\) and \(\gamma'''\),
respectively. The \(l\)-th derivative of \(\gamma\) with respect to \(s\) will be written as \(\gamma^{(l)}\)
where \(l=4,\ldots,2r\).

\par\medskip
\textbf{Acknowledgements:}
The author would like to thank Stefano Montaldo, Cezar Oniciuc and Andrea Ratto for many helpful discussions on the topic and the two referees for their many helpful suggestions
which helped to substantially improve the article.

\section{Proofs of the main results}
In this section we provide the proofs of the main results of this article,
which are Theorems \ref{theorem-relation-kt-3d} and \ref{theorem-r-harmonic-a}.

In the case of a space form the Riemann curvature tensor acquires the simple form
\begin{align*}
R(X,Y)Z=K(\langle Z,Y\rangle X-\langle Z,X\rangle Y),
\end{align*}
where \(X,Y,Z\) are vector fields and \(K\) the sectional curvature.
Then, the equation for polyharmonic curves \eqref{el-poly} simplifies as follows
\begin{align}
\label{polyharmonic-space-form}
\tau_r(\gamma)=\nabla_T^{2r-1}T+K\sum_{l=0}^{r-2}(-1)^l\big(\langle T,\nabla^l_TT\rangle\nabla_T^{2r-3-l}T
-\langle T,\nabla_T^{2r-3-l}T\rangle\nabla_T^lT\big)=0,
\end{align}
where \(T=\gamma'\) represents the unit tangent vector of the curve \(\gamma\).

Now, we assume that the target manifold \(N\) is a three-dimensional space form with a metric of constant 
curvature \(K\). Moreover, let \(\gamma\) be an arclength parametrized curve with tangent vector \(T\).
We choose the orthonormal Frenet frame \(\{T,F_2,F_3\}\) which satisfies the Frenet equations
\begin{align}
\label{frenet-3d}
\nabla_TT=kF_2,\qquad \nabla_TF_2=-kT+\tau F_3,\qquad \nabla_TF_3=-\tau F_2,
\end{align}
where \(k\) represents the geodesic curvature of the curve \(\gamma\)
and \(\tau\) its torsion. Again, note that we use the symbol \(\tau\) to denote both the torsion of a curve
and the tension field of a map, see \eqref{tension-field}.

\begin{Lem}
Let \(\gamma\colon I\to N\) be a curve parametrized by arclength together
with its Frenet frame \(\{T,F_2,F_3\}\). Furthermore, we assume that \(k\)
and \(\tau\) are constant. Then the following formulas hold
\begin{align}
\label{iterated-frenet-3d}
\nabla^{2l}_TT=&(-1)^lk^2(k^2+\tau^2)^{l-1}T+k\tau(-1)^{l+1}(k^2+\tau^2)^{l-1}F_3,\\
\nonumber\nabla^{2l+1}_TT=&(-1)^lk(k^2+\tau^2)^{l}F_2,
\end{align}
where \(l\geq 1\).
\end{Lem}

\begin{proof}
Using the Frenet-equations \eqref{frenet-3d} we find
\begin{align*}
\nabla^2_TT=&-k^2T+k\tau F_3,\\
\nabla^3_TT=&-k(k^2+\tau^2)F_2,\\
\nabla^4_TT=&k^2(k^2+\tau^2)T-k\tau(k^2+\tau^2)F_3,\\
\nabla^5_TT=&k(k^2+\tau^2)^2F_2,\\
\nabla^6_TT=&-k^2(k^2+\tau^2)^2T+k\tau(k^2+\tau^2)^2F_3,\\
\nabla^7_TT=&-k(k^2+\tau^2)^3F_2,\\
\nabla^8_TT=&k^2(k^2+\tau^2)^3T-k\tau(k^2+\tau^2)^3F_3.
\end{align*}
The claim then follows by iteration.
\end{proof}

\begin{proof}[Proof of Theorem \ref{theorem-relation-kt-3d}]
First, we assume that \(r\) is even. 
Using \eqref{iterated-frenet-3d} in \eqref{polyharmonic-space-form} we find
\begin{align*}
\sum_{l=0}^{r-2}(-1)^l\langle T,\nabla^l_TT\rangle\nabla_T^{2r-3-l}T
&=\sum_{l=0}^{\frac{r-2}{2}}\langle T,\nabla^{2l}_TT\rangle\nabla^{2r-3-2l}_TT \\
&=k(k^2+\tau^2)^{r-2}F_2+(\frac{r}{2}-1)k^3(k^2+\tau^2)^{r-3}F_2, \\
\sum_{l=0}^{r-2}(-1)^l\langle T,\nabla_T^{2r-3-l}T\rangle\nabla_T^lT
&=-\sum_{l=0}^{\frac{r-4}{2}}\langle T,\nabla^{2r-2l-4}_TT\rangle\nabla^{2l+1}_TT \\
&=(1-\frac{r}{2})k^3(k^2+\tau^2)^{r-3}F_2.
\end{align*}

Hence, we may conclude that the equation for a \(r\)-harmonic curve \eqref{polyharmonic-space-form} implies
\begin{align*}
0=k(k^2+\tau^2)^{r-3}\big(-(k^2+\tau^2)^2+K(k^2+\tau^2+(r-2)k^2)\big)F_2
\end{align*}
yielding the claim in the case that \(r\) is even.
The case of \(r\) being odd can be treated similarly.
\end{proof}

\begin{proof}[Proof of Corollary \ref{cor-helices}]
An arbitrary helix in \(\s^3\) can be parametrized by
\begin{align}
\label{helix-s3-general}
\gamma(s)=\big(\cos\alpha\cos(as),\cos\alpha\sin(as),\sin\alpha\cos(bs),\sin\alpha\sin(bs)\big),
\end{align}
where \(\alpha\in(0,\pi/2)\) and \(a,b\) are positive real numbers.
We require that \(a^2\cos^2\alpha+b^2\sin^2\alpha=1\)
which ensures that \(|\gamma'|^2=1\).
We choose \(a>b\) and find that the geodesic curvature \(k\) and the torsion \(\tau\)
of the curve \eqref{helix-s3-general} are given by
\begin{align}
\label{helix-s3-curvature-torsion}
k=\sqrt{(a^2-1)(1-b^2)},\qquad \tau=ab.
\end{align}
For more details on the geometry of helices in \(\s^3\) we refer to \cite[Section 3]{MR979640}.
Inserting \eqref{helix-s3-curvature-torsion} into the general formula
relating torsion and curvature of a polyharmonic curve \eqref{relation-kt-3d}
we find after a straightforward calculation
\begin{align}
\label{system-a-b}
(a^2+b^2-1)^2=(r-1)(a^2+b^2-a^2b^2-1)+a^2b^2.
\end{align}
The system \eqref{system-a-b} describes a particular conic section
and one can expect that there exist many positive \(a,b\in\R\) that
satisfy \eqref{system-a-b}. 
In order to find one particular solution we set
\begin{align*}
a=R\cos(t),~~b=R\sin(t),\qquad R>1,~~ t\in(0,\epsilon).
\end{align*}
Inserting into \eqref{system-a-b} then gives
\begin{align*}
\sin(2t)=\frac{2}{R^2}\sqrt{\frac{(R^2-1)(r-R^2)}{r-2}},
\end{align*}
where we now assume that \(r\geq 3\) and \(R^2<r\).
For \(R^2\) being close to \(r\) the right-hand side of this equation 
is close to zero such that there exists a \(t\in(0,\epsilon)\)
solving
\begin{align*}
t(R):=\frac{1}{2}\arcsin\big(\frac{2}{R^2}\sqrt{\frac{(R^2-1)(r-R^2)}{r-2}}\big).
\end{align*}
It is straightforward to check that
\begin{align*}
\lim_{R^2\to r}R\cos(t(R))=\sqrt{r},\qquad \lim_{R^2\to r}R\sin(t(R))=0
\end{align*}
and hence for \(R^2\) close enough to \(r\) we get
\(a=R\cos(t)>1\) and \(b=R\sin(t)<1\) that solve the system \eqref{system-a-b}.
Finally, as \(a>1,b<1\), 
one can always find a number \(\alpha\in(0,\pi/2)\) such that
the constraint \(a^2\cos^2\alpha+b^2\sin^2\alpha=1\)
is satisfied.
\end{proof}

In the following we will provide the proof of Theorem \ref{theorem-r-harmonic-a}.
We consider an ansatz of the form
\begin{align*}
\gamma(s)=\cos(as)e_1+\sin(as)e_2+e_3,
\end{align*}
where \(e_i,i=1,2,3\) are mutually perpendicular, \(|e_1|^2=|e_2|^2=\alpha^2,|e_1|^2+|e_3|^2=1\) and \(a\in\R\).
Note that we do not require the curve to be parametrized by arclength at this point.

We will use the 
relation between the \(r\)-energy \eqref{energy-r-harmonic} and the associated Lagrangian
\begin{align*}
E_r(\gamma)=\int_I\mathcal{L}_r^{\s^n}ds=\int_I|\nabla^{r-1}_TT|^2ds.
\end{align*}
Moreover, we will make use of the inclusion map \(\iota\colon\s^n\to\R^{n+1}\)
and also exploit the special structure of the Levi-Civita connection \(\nabla\) on the sphere
\begin{align*}
d\iota(\nabla_TX)=X'+\langle X,\gamma'\rangle\gamma,
\end{align*}
where \(X\) is a vector field on \(\s^n\subset\R^{n+1}\).

Using these facts we give an expression of the Lagrangian associated with the \(r\)-energy for arbitrary values of \(r\).

\begin{Lem}
Consider the curve \(\gamma\colon I\to\s^n\) given by
\begin{align}
\label{ansatz-kharmonic-a}
\gamma(s)=\cos(as)e_1+\sin(as)e_2+e_3,
\end{align}
where \(e_i,i=1,2,3\) are mutually perpendicular, \(|e_1|^2=|e_2|^2=\alpha^2,|e_1|^2+|e_3|^2=1\) and \(a\in\R\).

Then, for any \(l\in\N\), the following formulas hold
\begin{align}
\label{r-tension-ansatz-kharmonic-a}
d\iota(\nabla^{2l}_TT)&=(-1)^la^{2l}(1-\alpha^2)^l\gamma',\\
\nonumber d\iota(\nabla^{2l+1}_TT)&=(-1)^la^{2l}(1-\alpha^2)^l d\iota(\nabla_TT),
\end{align}
where \(\iota\colon\s^n\to\R^{n+1}\) denotes the embedding of \(\s^n\) into \(\R^{n+1}\). 
\end{Lem}
\begin{proof}
Using the above ansatz a direct calculation shows
\begin{align*}
d\iota(\nabla_TT)&=\gamma''+a^2\alpha^2\gamma,\\
d\iota(\nabla^2_TT)&=\gamma'''+a^2\alpha^2\gamma' \\
&=-a^2(1-\alpha^2)\gamma',\\
d\iota(\nabla^3_TT)&=\gamma^{(4)}
+\langle\gamma''',\gamma'\rangle\gamma+a^2\alpha^2\gamma''+a^4\alpha^4\gamma ,\\
&=-a^2(1-\alpha^2)(\gamma''+a^2\alpha^2\gamma)\\
&=-a^2(1-\alpha^2)d\iota(\nabla_TT),\\
d\iota(\nabla^4_TT)&=
a^4(1-\alpha^2)^2\gamma'
,\\
d\iota(\nabla^5_TT)&=a^4(1-\alpha^2)^2d\iota(\nabla_TT).
\end{align*}
The claim now easily follows by induction.
\end{proof}

Using \eqref{r-tension-ansatz-kharmonic-a} we are now able to give a recursion formula
for the \(r\)-energy of a curve of the form \eqref{ansatz-kharmonic-a}.
Similar recursion formulas for the inclusion map \(\iota:\s^n\hookrightarrow\s^{n+1}\) were
obtained in \cite[Lemma 3.8]{MR3711937} in the study of proper \(r\)-harmonic submanifolds
and in \cite[Proposition 2.10]{MR4106647} in the study of proper \(ES-r\)-harmonic submanifolds.

\begin{Lem}
Consider the curve \(\gamma\colon I\to\s^n\) given by
\begin{align*}
\gamma(s)=\cos(as)e_1+\sin(as)e_2+e_3,
\end{align*}
where \(e_i,i=1,2,3\) are mutually perpendicular, \(|e_1|^2=|e_2|^2=\alpha^2,|e_1|^2+|e_3|^2=1\) and \(a\in\R\).
Then the Lagrangian associated with the \(r\)-energy of \(\gamma\) is given by
\begin{align}
\label{recursion-r-energy}
\mathcal{L}^{\s^n}_{r}(\alpha)=a^{2r}\alpha^2(1-\alpha^2)^{r-1},
\end{align}
where \(r\geq 1\).
\end{Lem}
\begin{proof}
Using \eqref{r-tension-ansatz-kharmonic-a} a direct calculation yields
\begin{align*}
|d\iota(\nabla^{2l}_TT)|^2&=a^{4l}(1-\alpha^2)^{2l}|\gamma'|^2\\
&=\alpha^2a^{4l+2}(1-\alpha^2)^{2l},\\
|d\iota(\nabla^{2l+1}_TT)|^2&
=a^{4l}(1-\alpha^2)^{2l}|d\iota(\nabla_TT)|^2\\
&=\alpha^2 a^{4l+4}(1-\alpha^2)^{2l+1}
\end{align*}
completing the proof.
\end{proof}

\begin{proof}[Proof of Theorem \ref{theorem-r-harmonic-a}]
By making use of the ansatz \eqref{ansatz-kharmonic-a} we have reduced the general \(r\)-energy 
for a curve \eqref{energy-poly-curve} to the simpler form \eqref{recursion-r-energy}.

Thus, we now consider the variation of the Lagrangian associated to \eqref{recursion-r-energy} with respect to \(\alpha\) and find
\begin{align*}
\frac{d}{d\alpha}\mathcal{L}_r^{\s^n}(\alpha)=2a^{2r}\alpha(1-\alpha^2)^{r-2}(1-\alpha^2r).
\end{align*}
The right hand side of this equation vanishes for \(\alpha=1\) (assuming that \(r\geq 3\)) corresponding to the 
case of \(\gamma\) being a geodesic but also for \(1=\alpha^2r\) which
is the case we are interested in.

At this point, we impose that the curve \(\gamma\) is parametrized with respect to arclength.
Employing the ansatz \eqref{ansatz-kharmonic-a} we find \(|\gamma'|^2=a^2\alpha^2\).
Combining the arclength constraint \(a^2\alpha^2=1\) with the equality obtained from the first variation formula \(1=\alpha^2r\)
we obtain \(a^2=r\) completing the proof.
\end{proof}

\begin{Bem}
As the curve \eqref{ansatz-kharmonic-a} is a plane curve we could choose \(\s^2\) instead of \(\s^n\)
in the proof of Theorem \ref{theorem-r-harmonic-a}.

In \cite[Remark 2.9]{MR3990379} it is noted that the composition of a biharmonic and a totally geodesic map is biharmonic
and one should expect that the same also holds true for \(r\)-harmonic maps.
Consequently, it should be possible to prove Theorem \ref{theorem-r-harmonic-a} by only working on \(\s^2\)
and then using that \(\s^2\) is totally geodesic in \(\s^n\) for \(n\geq 3\).
\end{Bem}

Let us also make a short comment on the stability of the family of curves \eqref{solution-kharmonic-a} given in Remark \ref{remark-general-solution}.

\begin{Prop}
Let \(\gamma\colon I\to\s^n\) be an \(r\)-harmonic curve parametrized by arclength of the form \eqref{solution-kharmonic-a}.
Then the second variation of the \(r\)-energy evaluated at a critical point is given by
\begin{align}
\label{solution-kharmonic-a-sv}
\frac{d^2}{d\alpha^2}E_r(\alpha)\big|_{\alpha^2r=1}=-4|I|r^{r}\big(\frac{r-1}{r}\big)^{r-2}<0,
\end{align}
where \(|I|=\int_Ids\).
\end{Prop}
\begin{proof}
By a direct calculation we find
\begin{align*}
\frac{d^2}{d\alpha^2}E_r(\alpha)=&2|I|a^{2r}(1-\alpha^2)^{r-3}(1-\alpha^2r)(1+3\alpha^2-2\alpha^2r)
-4|I|a^{2r}\alpha^2r(1-\alpha^2)^{r-2}.
\end{align*}
Evaluating this formula at the critical point \(\alpha^2r=1,a^2=r\) completes the proof.
\end{proof}

\bibliographystyle{plain}
\bibliography{mybib}
\end{document}